\documentclass{amsart}
\usepackage{}
\usepackage{cases}
\usepackage{amsmath}
\usepackage{amssymb}
\usepackage{amsfonts}
\newtheorem{theorem}{Theorem}[section]
\newtheorem{lemma}[theorem]{Lemma}
\newtheorem{corollary}[theorem]{Corollary}
\newtheorem{proposition}[theorem]{Proposition}

\theoremstyle{definition}

\newtheorem{remark}[theorem]{Remark}
\numberwithin{equation}{section}

\begin{document}
\title[$L^p$-Liouville theorems on metric measure spaces]
{$L^p$-Liouville theorems on complete smooth\\ metric measure spaces}

\author{Jia-Yong Wu}
\address{Department of Mathematics, Shanghai Maritime University,
Haigang Avenue 1550, Shanghai 201306, P. R. China}
\email{jywu81@yahoo.com}

\thanks{This work was partially supported by NSFC
(11101267, 11271132) and the Innovation Program of Shanghai Municipal
Education Commission (13YZ087).}


\dedicatory{}

\date{\today}

\begin{abstract}
We study some function-theoretic properties on a complete smooth
metric measure space $(M,g,e^{-f}dv)$ with Bakry-\'{E}mery
Ricci curvature bounded from below. We derive a Moser's parabolic Harnack
inequality for the $f$-heat equation, which leads to upper
and lower Gaussian bounds on the $f$-heat kernel. We also
prove $L^p$-Liouville theorems in terms of the lower bound of
Bakry-\'{E}mery Ricci curvature and the bound of function $f$,
which generalize the classical Ricci curvature case and
the $N$-Bakry-\'{E}mery Ricci curvature case.
\end{abstract}
\maketitle

\section{Introduction and main results}
\subsection{Background}
Let $(M,g)$ be an $n$-dimensional complete Riemmannian manifold
and $f$ be a smooth function on $M$. We define a symmetric
diffusion operator $\Delta_f$ (or $f$-Laplacian), which is given by
\[
\Delta_f:=\Delta-\nabla f\cdot\nabla,
\]
where $\Delta$ and $\nabla$ are the Laplacian and covariant derivative
of the metric $g$. The $f$-Laplacian $\Delta_f$ is the infinitesimal
generator of the Dirichlet form
\[
\mathcal {E}(\phi_1,\phi_2)=\int_M (\nabla \phi_1, \nabla \phi_2)d\mu,
\,\,\,\forall \phi_1, \phi_2\in C_0^{\infty}(M),
\]
where $\mu$ is an invariant measure of $\Delta_f$ given by $d\mu=e^{-f}dv$,
and where $dv$ is the volume element induced by the metric $g$. It is
well-known that $\Delta_f$ is self-adjoint with respect to the weighted
measure $d\mu$. The triple $(M,g,e^{-f}dv)$ is customarily called
a complete smooth metric measure space. On this measure space,
we often consider the $f$-heat equation
\[
\left(\frac{\partial}{\partial t}-\Delta_f\right)u=0
\]
instead of the classical heat equation. If the function $u$
is independent of time $t$, then $u$ is a $f$-harmonic
function. In this paper, we denote by $H(x,y,t)$ the $f$-heat
kernel, that is, for each $x\in M$, $H(x,y,t)=u(y,t)$ is
the minimal solution of the $f$-heat equation with
$u(y,0)=\delta_x(y)$. Equivalently, $H(x,y,t)$ is the kernel
of the semigroup $P_t=e^{t\Delta_f}$ associated to the
Dirichlet form $\mathcal {E}(\phi,\phi)$.

On the smooth metric measure space $(M,g,e^{-f}dv)$, following Bakry and
\'{E}mery  \cite{[Ba]} and \cite{[BE]} (see also \cite{[LD]} and
\cite{[Lott1]}), we define the Bakry-\'{E}mery Ricci curvature
\[
Ric_f:=Ric+Hess(f),
\]
where $Ric$ denotes the Ricci curvature of the manifold and $Hess$
denotes the Hessian with respect to the manifold metric. The
Bakry-\'{E}mery Ricci curvature is a natural extension of the Ricci
curvature. If $f$ is constant, $Ric_f$ returns to the Ricci
curvature $Ric$. The Bakry-\'{E}mery Ricci curvature has been
extensively studied because it often shares interesting properties
with the ordinary Ricci curvature. For example, there exists an
interesting Bochner type identity
\[
\Delta_f|\nabla u|^2=2|Hess(u)|^2+2\langle\nabla u,\nabla\Delta_fu\rangle
+2Ric_f(\nabla u,\nabla u).
\]
This identity is parallel to the Bochner identity in the classical
Ricci curvature case, and plays
an important role in studying comparison theorems (see \cite{[WW]}).
For more extended results, the interested reader can consult
\cite{[BQ1]}, \cite{[BQ2]}, \cite{[CSW]}, \cite{[FLZ]}, \cite{[LD2]},
\cite{[Lott1]}, \cite{[Lo-Vi]}, \cite{[St1]}, \cite{[St2]},
\cite{[Wang]} and \cite{[Wu]}.

Also, the Bakry-\'{E}mery Ricci curvature has become an important object of
study in geometry analysis, in large part due to so-called
gradient Ricci solitons. Recall that a complete manifold $(M,g)$ is
a gradient Ricci soliton if for some function $f$ on $M$ and some
constant $\rho$ we have that
\[
Ric_f=\rho g.
\]
The soliton is called expanding, steady or shrinking if,
respectively, $\rho<0$, $\rho=0$
and $\rho>0$. Ricci solitons possess many interesting
geometric and topological properties. See, for example, \cite{[Cao1]},
\cite{[Cao2]} and \cite{[Perelman]} for nice explanations on
this subject.

Recently, there has been renewed interest in the Bakry-\'{E}mery Ricci
curvature and its modified version, the $N$-Bakry-\'{E}mery Ricci
curvature, defined by
\[
Ric^N_f:=Ric+Hess(f)-\frac{df\otimes df}{N},
\]
where $N$ is a positive constant. For example, Catino,
Mantegazza, Mazzieri and Rimoldi \cite{[CMMR]},
Petersen and Wylie \cite{[PeWy]}, and Pigola,
Rimoldi and Setti \cite{[PiRiSe]} established various
Liouville-type or rigidity theorems about these curvatures.
Prior to their works, X.-D. Li \cite{[LD]} studied a
$L^1$-Liouville theorem in case the $N$-Bakry-\'{E}mery
Ricci curvature is bounded below by a negative quadratic
polynomial of the distance function. That is an extension
of the classical $L^1$-Liouville theorem on Ricci curvatures,
proved by P. Li \cite{[Li0]}. However, as X.-D. Li said
in Subsection 8.6 of \cite{[LD]}, we cannot prove a
$L^1$-Liouville theorem if we only assume a lower bound of
the same kind on $Ric_f$. Indeed, we can not obtain a
Li-Yau type parabolic Harnack inequality under only this
curvature assumption. Here it is natural to pose the
following problem: What are the optimal geometric or
analytic conditions on the smooth metric measure space
in order that the Li-Yau parabolic Harnack inequality holds?

In the recent papers \cite{[MuWa],[MuWa2]}, Munteanu and Wang
partially answered to the above question. In particular, they
derived gradient estimates and Liouville properties for
positive $f$-harmonic functions under suitable growth
assumption on $f$. Their theorems take the form of Yau's
classical result on positive $f$-harmonic functions, but
the proof they adopt is new and quite different in spirit
from Yau's direct application of the Bochner formula
\cite{[Yau]}. Their approach essentially relies on the well-known
De Giorgi-Nash-Moser theory. This motivates our proof
of Theorem \ref{THEO} in this paper.

\subsection{Main results}\label{subb2.1}
The purpose of this paper is to further study geometric
inequalities for the $f$-heat equation and $L^p$-Liouville
theorems for $f$-harmonic functions on complete smooth
metric measure spaces. One contribution of this paper is to
provide suitable weighted curvature conditions which assure
the validity of various well-known geometric inequalities,
such as a local $f$-volume doubling property, a local
$f$-Neumann Poincar\'e inequality and a local $f$-mean value
inequality, etc.. Another contribution of this paper is that
we used those geometric inequalities to prove new
$L^p$-Liouville theorems on complete smooth metric measure
spaces.

This paper can be divided into two parts. In the first part,
borrowing the idea of Munteanu and Wang \cite{[MuWa],[MuWa2]},
we will derive some geometric inequalities, such as parabolic
Harnack inequalities, H\"{o}lder continuity estimates and
$f$-heat kernel estimates on complete smooth metric measure
spaces. We first present a parabolic Harnack inequality
on complete smooth metric measure spaces.
\begin{theorem}\label{THEO}
Let $(M,g,e^{-f}dv)$ be an $n$-dimensional complete noncompact
smooth metric measure space. If $Ric_f\geq-(n-1)K$ and
$|f|(x)\leq A$ for some nonnegative constants $K$ and $A$,
then there exist a constant $c(n,A)$ such that, for any
$0<R\leq\infty$ and ball $B_o(r)$, $o\in M$, $0<r<R$ and for
any smooth positive solution $u$ of the $f$-heat equation in the
cylinder $Q=B_o(r)\times (s-r^2,s)$, we have
\[
\sup_{Q_{-}}\{u\}\leq e^{c(n,A)(1+Kr^2)}\cdot \inf_{Q_{+}}\{u\},
\]
where $Q_{-}:=B_o(\frac 12r)\times (s-\frac 34r^2,s-\frac 12r^2)$
and $Q_{+}:=B_o(\frac 12r)\times (s-\frac 14r^2,s)$.
\end{theorem}
The sketch of the proof of Theorem \ref{THEO} will be given in Section
\ref{sec2a}. The proof follows by the Moser iteration technique
\cite{[Moser]}, which involves a local Sobolev inequality on a
smooth metric measure space. Munteanu and Wang \cite{[MuWa2]}
used a similar technique to derive an elliptic Harnack inequality
for $f$-harmonic functions. When the metric measure space is a
Riemannian manifold, that is, the function $f$ is constant, this
result was obtained independently by Saloff-Coste \cite{[Saloff]}
and Grigor'yan \cite{[Grig]}.

A standard consequence of Theorem \ref{THEO} is a strong
Liouville theorem for any $f$-harmonic function (see Corollary
\ref{holliou}). Theorem \ref{THEO} also implies two-sided
$f$-heat kernel bounds. This result is essentially analogous
to the case of heat equation on Riemannian manifolds in
\cite{[Saloff2]} (see also \cite{[Grig3]}).
\begin{theorem}\label{MainAn}
Let $(M,g,e^{-f}dv)$ be an $n$-dimensional complete noncompact
smooth metric measure space. If $Ric_f\geq-(n-1)K$ and
$|f|(x)\leq A$ on the ball $B_o(2R)$ for some nonnegative
constants $K$ and $A$, then there exist positive constants
$c_i$, $i=5,6,7,8$, depending only on $n$ and $A$ such that
\[
\frac{e^{-c_6(1+Kt)}}{V_f(B_x(\sqrt{t}))}\exp\left(-c_5\frac{d^2(x,y)}{t}\right)
\leq H(x,y,t)\leq \frac{e^{c_8(1+Kt)}}{V_f(B_x(\sqrt{t}))}
\exp\left(-c_7\frac{d^2(x,y)}{t}\right)
\]
for any $x,y\in B_o(R/2)$ and $0<t<R^2/4$, where $V_f(B_x(\sqrt{t}))$ denotes
the $f$-volume of the ball $B_x(\sqrt{t})$ with respect to $e^{-f}dv$.
\end{theorem}
\begin{remark}
Theorem \ref{MainAn} gives an accurate description of the
coefficients of two-sided $f$-heat kernel bounds.
It will be crucial in the proof of Theorem \ref{Main3}.
\end{remark}

The proof strategy of Theorem \ref{MainAn}
is different from the classical Li-Yau trick \cite{[Li0]}.
In \cite{[Li0]}, two-sided Gaussian bounds on the heat
kernel are obtained by the Li-Yau gradient estimate.
However, in our case it seems to be impossible to adopt
Li-Yau gradient estimate method directly in order to
derive upper and lower bounds on the $f$-heat kernel on
complete smooth metric measure spaces. In our approach,
Gaussian bounds on the $f$-heat kernel rely on the Moser's
parabolic Harnack inequality and the integral estimate
of the $f$-heat kernel due to Davies \cite{[Davies]},
thus our arguments are similar to the ones of Saloff-Coste
\cite{[Saloff],[Saloff1],[Saloff2]} and Grigor'yan
\cite{[Grig]}. Please see Section \ref{sec2c} for a
detailed discussion.

\vspace{0.5em}

In the second part of this paper, we will investigate
various $L^p$-Liouville theorems for $f$-harmonic functions
on complete noncompact metric measure space $(M,g,e^{-f}dv)$
under different assumptions on $Ric_f$ and $f$.

\vspace{0.5em}

We first start recalling a $L^p$-Liouville theorem for positive
$f$-subharmonic functions when $1<p<\infty$, which extends the
result in the classical case in \cite{[Yau2]}.
This was originally proved
in \cite{[PiRiSe05]}; see also \cite{[PiRiSe]}.
\begin{theorem}[Pigola, Rigoli and Setti \cite{[PiRiSe05]}]\label{Main1}
Let $(M,g,e^{-f}dv)$ be an $n$-dimensional complete smooth metric measure
space. For any $1<p<\infty$, there does not exist any nonconstant,
nonnegative, $L^p(\mu)$-integrable $f$-subharmonic function.
\end{theorem}

We now deal with the $L^p$-Liouville theorem in case of $0<p<1$.
In this case we obtain an analogous result to that obtained
by Li and Schoen in \cite{[Li-Sch]}. See Subsection \ref{sub6.2}
for a detailed discussion.

\begin{theorem}\label{Main2}
Let $(M,g,e^{-f}dv)$ be an $n$-dimensional complete noncompact
smooth metric measure space. Assume that $f$ is bounded, and
there exists a constant $\delta(n)>0$ depending only on $n$,
such that, for some point $o\in M$, the
Bakry-\'{E}mery Ricci curvature satisfies
\[
Ric_f\geq-\delta(n)r^{-2}(x),
\]
whenever the distance from $o$ to $x$, $r(x)$, is sufficiently
large. Then any nonnegative $L^p(\mu)$-integrable ($0<p<1$)
$f$-subharmonic function must be identically zero.
\end{theorem}

\vspace{0.5em}

Finally, motivated by the P. Li's work \cite{[Li0]} and
X.-D. Li's generalization \cite{[LD]}, we obtain a new
$L^1$-Liouville theorem on smooth metric measure spaces.
\begin{theorem}\label{Main3}
Let $(M,g,e^{-f}dv)$ be an $n$-dimensional complete noncompact
smooth metric measure space. Assume that $f$ is bounded, and
there exists a constant $C>0$, such that, for some point $o\in M$,
the Bakry-\'{E}mery Ricci curvature satisfies
\[
Ric_f\geq-C(1+r^2(x)),
\]
where $r(x)$ denotes the distance from $o$ to $x$. Then any
nonnegative $L^1(\mu)$-integrable $f$-subharmonic function
must be identically constant.
\end{theorem}

Theorem \ref{Main3} partially answers to a question posed by
X.-D. Li (see Subsection 8.6 in \cite{[LD]}). Its proof is
similar to the arguments of \cite{[Li0]}, where a critical
step is the usage of the upper Gaussian bound on the
$f$-heat kernel proved in Theorem \ref{MainAn}. A detailed
discussion shall be carried out in Subsection \ref{sub6.3}.

\begin{remark}
We remark that the absolute value of a $f$-harmonic function
is a nonnegative $f$-subharmonic. Therefore we can conclude that
a complete metric measure space does not admit any nonconstant
$L^p(\mu)$-integrable $f$-harmonic function under the same
hypotheses of Theorems \ref{Main1}, \ref{Main2} and
\ref{Main3}, respectively.
\end{remark}

\begin{remark}
As many recent authors said in \cite{[DuSu]}, \cite{[SuZh]} and
\cite{[Wu2]}, if the condition on $f$ bounded is replaced by
$|\nabla f|$ bounded, then similar results to Theorems \ref{Main2}
and \ref{Main3} can be immediately obtained by modifying the
arguments of \cite{[LD]}. Indeed, the conditions $Ric_f\geq-(n-1)K$
and $|\nabla f|\leq a$ imply that
\[
Ric^N_f\geq-(n-1)\left(K+\frac{a^2}{N(n-1)}\right).
\]
\end{remark}

The rest of this paper is organized as follows. In Section
\ref{sec2a}, we present a local $f$-volume doubling property,
a local $f$-Neumann Poincar\'e inequality and a local Sobolev
inequality on complete smooth metric measure spaces. After that,
following the arguments of Saloff-Coste \cite{[Saloff]} or
Grigor'yan \cite{[Grig]}, we establish a Moser's version of
parabolic Harnack inequality. In Section \ref{sec2b},
using the parabolic Harnack inequality, we obtain a
H\"{o}lder continuity estimate for the $f$-heat equation,
which implies a strong Liouville theorem.
In Section \ref{sec2c}, we prove two-sided Gaussian bounds on
the $f$-heat kernel on complete smooth metric measure spaces.
In Section \ref{sec3}, we derive a $f$-mean value inequality on
complete smooth metric measure spaces, which is similar to the
case of harmonic functions on a manifold, obtained by Li and
Schoen \cite{[Li-Sch]}. In Section \ref{sec4}, we establish
$L^p$-Liouville theorems on complete smooth metric measure
spaces by following the ideas in \cite{[Li0]} and \cite{[Li-Sch]}.

\section{Poincar\'e, Sobolev and Harnack inequalities}
\label{sec2a}
Let $\Delta_f=\Delta-\nabla f\cdot\nabla$ be the $f$-Laplacian
on a complete smooth metric measure space $d\mu=e^{-f}dv$ on a
complete Riemannian manifold. For a set $\Omega$, we will
denote by $V(\Omega)$ the volume, and by $V_f(\Omega)$ the
$f$-volume of $\Omega$. Throughout this section, we will assume
\[
Ric_f\geq -(n-1)K\quad
\mathrm{and}\quad |f|(x)\leq A
\]
for some nonnegative constants $K$ and $A$. Under these assumptions,
we have the validity of the $f$-Laplacian and $f$-volume
comparison results.
\begin{lemma}[Wei and Wylie \cite{[WW]}]\label{comp}
Let $(M,g,e^{-f}dv)$ be an $n$-dimensional complete noncompact
smooth metric measure space. If $Ric_f\geq-(n-1)K$ and
$|f|(x)\leq A$ for some nonnegative constants $K$ and $A$, then
along any minimizing geodesic starting from $x\in M$ we have

\[
\Delta_fr(x,y)\leq(n-1+4A)\sqrt{K}\coth \sqrt{K}r
\]
for any $0<r<R$, where $r(x,y):=d(x,y)$ is the distance function.
Hence along any minimizing geodesic starting from
$x\in M$ we have
\begin{equation}\label{volcomp}
\frac{V_f(B_x(r_2))}{V_f(B_x(r_1))}
\leq\frac{V^{n+4A}_K(r_2)}{V^{n+4A}_K(r_1)}
\end{equation}
for any $0<r_1<r_2<R$. Here $V^{n+4A}_K(r)$ is the volume
of the radius $r$-ball in the model space $M^{n+4A}_K$,
the simply connected model space of dimension $n+4A$ with
constant curvature $K$.
\end{lemma}

From \eqref{volcomp}, we easily deduce that
\begin{equation}\label{voldop}
V_f(B_x(2r))\leq 2^{n+4A}e^{C(n,A)\sqrt{K}r}\cdot V_f(B_x(r))
\end{equation}
for any $0<r<R$. This inequality implies that the local
$f$-volume doubling property holds. This property will
play a crucial role in our paper. We say that a complete smooth metric
measure space $(M,g,e^{-f}dv)$ admits a \emph{local $f$-volume
doubling property} if for any fixed $0<R<\infty$, there exists
a constant $C(R)$ such that
\[
V_f(B_x(2r))\leq C(R)\cdot V_f(B_x(r))
\]
for any $0<r<R$ and $x\in M$. Note that, when $K=0$, the above
inequality holds with $R=+\infty$, and it called the \emph{global
$f$-volume doubling property}.

By Lemma \ref{comp}, following the Buser's proof \cite{[Bus]} or
the Saloff-Coste's alternate proof (Theorem 5.6.5  in
\cite{[Saloff2]}), we can easily get a local Neumann Poincar\'e
inequality in the setting of smooth metric measure spaces.
\begin{lemma}\label{NeuPoin}
Let $(M,g,e^{-f}dv)$ be an $n$-dimensional complete noncompact
smooth metric measure space, and denote by $r(x)$ the distance
function from a fixed origin $o\in M$. If $Ric_f\geq-(n-1)K$ and
$|f|(x)\leq A$ for some nonnegative constants $K$ and $A$,
then
\begin{equation}\label{Nepoinineq}
\int_{B_o(r)}|\varphi-\varphi_{B_o(r)}|^2e^{-f}dv\leq
e^{c_1(1+\sqrt{K}r)}\cdot r^2\int_{B_o(r)}|\nabla \varphi|^2e^{-f}dv
\end{equation}
for any $x\in M$ such that $0<r(x)<R$ and $\varphi\in C^\infty(B_o(r))$, where
$\varphi_{B_o(r)}:=V_f^{-1}(B_o(r))\int_{B_o(r)}\varphi e^{-f}dv$. The
constant $c_1$ depends only on the dimension $n$ and $A$.
\end{lemma}

\begin{remark}
Inequality \eqref{Nepoinineq} implies that a \emph{local $f$-Neumann
Poincar\'e inequality} holds. In \cite{[MuWa]}, Munteanu and Wang
proved a $f$-Neumann Poincar\'e inequality when $Ric_f\geq 0$. In
\cite{[MuWa2]}, they only obtained the validity of a $f$-Neumann
Poincar\'e inequality uniformly at small scales. But in our case,
the $f$-Poincar\'e inequality can hold on balls of \emph{any} radius
due to a stronger assumption on $f$, which is a crucial
step on the proof of $L^1$-Liouville theorem. Because in the
course of proof of $L^1$-Liouville result, we need to let the
radius of balls tend to infinity. Also note that when $f$ is
constant, \eqref{Nepoinineq} was obtained by Saloff-Coste
(see (6) in \cite{[Saloff1]} or Theorem 5.6.5 in \cite{[Saloff2]}).
\end{remark}

Combining Lemma \ref{comp}, Lemma \ref{NeuPoin} and
the argument in \cite{[Saloff]}, we have a local Sobolev
inequality, which is one of the key technical points needed
to apply Moser's iterative technique to derive parabolic
Harnack inequalities for the $f$-heat equation.
\begin{lemma}\label{NeuSob}
Let $(M,g,e^{-f}dv)$ be an $n$-dimensional complete noncompact
smooth metric measure space. If $Ric_f\geq-(n-1)K$ and
$|f|(x)\leq A$ for some nonnegative constants $K$ and $A$, then
for any constant $p>2$, there exists a constant $c_2$, depending
on $n$ and $A$ such that
\[
\left(\int_{B_o(r)}|\varphi|^{\frac{2p}{p-2}}e^{-f}dv
\right)^{\frac{p-2}{p}}\leq \frac{e^{c_2(1+\sqrt{K}r)}\cdot r^2}{V_f(B_o(r))^{\frac 2p}}
\int_{B_o(r)}(|\nabla \varphi|^2+r^{-2}|\varphi|^2)e^{-f}dv
\]
for any $x\in M$ such that $0<r(x)<R$ and $\varphi\in C^\infty(B_o(r))$.
\end{lemma}

\begin{proof}[Sketch proof of Lemma \ref{NeuSob}]
The proof is nearly the same as that of Theorem 2.1 in \cite{[Saloff]}
or Theorem 3.1 in \cite{[Saloff1]} except for our discussion with
respect to the weighted measure $e^{-f}dv$. When $f$ is constant,
this result was confirmed by Saloff-Coste \cite{[Saloff]}
(see also Theorem 3.1 in \cite{[Saloff1]}). We refer the
reader to these papers for a nice proof.
\end{proof}

\begin{remark}
In Lemma \ref{NeuSob}, the local Sobolev inequality is
different from Munteanu-Wang's Neumann Sobolev inequality
(Lemma 3.3 in \cite{[MuWa2]}). Here we mainly follow the
arguments of Saloff-Coste \cite{[Saloff]} to derive the
local Sobolev inequality, whereas Munteanu and Wang proved
their local Neumann Sobolev inequality adapting the same
arguments as in \cite{[HaKo]}. Note also that, while
Munteanu and Wang \cite{[MuWa2]} established the local Neumann
Sobolev inequality only on the \emph{unit} balls due to a weaker
hypothesis on the oscillation of $f$ on unit balls, our
local Sobolev inequality holds on balls of \emph{any} radius,
due to a stronger assumption on $f$.
\end{remark}

We shall now present a result concerning the Harnack inequality
for the $f$-heat equation, which is very much similar to the case
when $f$ is constant, obtained by Saloff-Coste \cite{[Saloff]}
and Grigor'yan \cite{[Grig]}.
\begin{theorem}\label{PoinDouHarn}
Let $(M,g,e^{-f}dv)$ be an $n$-dimensional complete noncompact
smooth metric measure space. Fix $0<R\leq\infty$. Assume that
\eqref{voldop} and \eqref{Nepoinineq} are satisfied up to this
$R$. Then there exist constants $c_3$ depending on $n$ and $A$
such that, for any ball $B_o(r)$, $o\in M$, $0<r<R$ and for any
smooth positive solution $u$ of the $f$-heat equation in the
cylinder $Q=B_o(r)\times (s-r^2,s)$, we have
\[
\sup_{Q_{-}}\{u\}\leq e^{c_3(1+Kr^2)}\cdot\inf_{Q_{+}}\{u\},
\]
where $Q_{-}:=B_o(\frac 12r)\times (s-\frac 34r^2,s-\frac 12r^2)$
and $Q_{+}:=B_o(\frac 12r)\times (s-\frac 14r^2,s)$.
\end{theorem}
\begin{proof}[Sketch proof of Theorem \ref{PoinDouHarn}]
The proof is the weighted case of the arguments of \cite{[Saloff]}
or \cite{[Saloff1]}. Indeed this result follows from the
standard Moser's technique. Since the conditions of Theorem
\ref{PoinDouHarn} imply a family of local Sobolev inequalities
due to Lemma \ref{NeuSob}, combining the local volume
doubling property, it is enough to run the Moser's iteration
procedure to prove Theorem \ref{PoinDouHarn}, as explained in
\cite{[Saloff]} or \cite{[Saloff1]}.
\end{proof}

Combining Lemmas \ref{comp}, \ref{NeuPoin} and Theorem
\ref{PoinDouHarn}, we immediately have that:
\begin{corollary}\label{coro1a}
Let $(M,g,e^{-f}dv)$ be an $n$-dimensional complete noncompact
smooth metric measure space. If $Ric_f\geq-(n-1)K$ and
$|f|(x)\leq A$ for some nonnegative constants $K$ and $A$,
then there exist a constant $c(n,A)$ such that for any
ball $B_o(r)$, $o\in M$, $0<r<R$ and for any smooth positive
solution $u$ of the $f$-heat equation in the
cylinder $Q=B_o(r)\times (s-r^2,s)$, we have
\[
\sup_{Q_{-}}\{u\}\leq e^{c(n,A)(1+Kr^2)}\cdot \inf_{Q_{+}}\{u\},
\]
where $Q_{-}:=B_o(\frac 12r)\times (s-\frac 34r^2,s-\frac 12r^2)$
and $Q_{+}:=B_o(\frac 12r)\times (s-\frac 14r^2,s)$.
\end{corollary}

\begin{remark}\label{rem28}
In \cite{[Saloff3]} and \cite{[Grig3]}, Saloff-Coste
and Grigor'yan have confirmed that the conjunction of the
$f$-volume doubling property and the $f$-Neumann Poincar\'e
inequality is equivalent to a parabolic Harnack inequality for
the $f$-heat equation. Our novel feature in this section is that
we take into account suitable weighted curvature condition which
implies the validity of these inequalities.
\end{remark}

\section{Liouville theorem}\label{sec2b}
In this section, we will apply the parabolic Harnack inequality
to obtain a quantitative H\"{o}lder continuity estimate for a
solution to the $f$-heat equation, and hence derive a strong
Liouville property under some suitable assumptions
on $Ric_f$ and $f$.

First, we give the H\"{o}lder continuity estimate for any
solution of the $f$-heat equation. When $f$ is constant this
was established in Theorem 5.4.7 of \cite{[Saloff2]}.
\begin{theorem}\label{MainHl}
Under the same assumptions of Theorem \ref{PoinDouHarn},
there exist $\theta\in(0,1)$, $\alpha\in(0,1)$ and
$A_\kappa=4\theta^{-1}(1-\kappa)^{-\alpha}>1$,
$\kappa\in(0,1)$, such that any solution $u$ of the
$f$-heat equation in $Q=B_o(r)\times (s-r^2,s)$, satisfies
\begin{equation*}
\sup_{(y,t),(y',t')\in Q_\kappa}
\left\{\frac{|u(y,t)-u(y',t')|}{[|t-t'|^{1/2}+d(y,y')]^\alpha}\right\}
\leq \frac{A_\kappa}{r^\alpha}\sup_Q\{|u|\},
\end{equation*}
where $Q_\kappa:=B_o(\kappa r)\times (s-\kappa r^2,s)$.
\end{theorem}
\begin{proof}
The proof is nearly the same as in \cite{[Moser]}
(see also \cite{[Saloff2]}) which uses the parabolic Harnack inequality.
For the reader's convenience, we include a detailed proof of this result.
For any non-negative solution $v$ of the $f$-heat equation in $Q$, by
Theorem \ref{PoinDouHarn}, we have
\begin{equation}\label{Harineq}
\frac{1}{\bar{V}_f(Q_-)}\int_{Q_{-}}vd\bar{\mu}\leq\max_{Q_{-}}\{v\}
\leq e^{c(n,A)(1+Kr^2)}\min_{Q_{+}}\{v\},
\end{equation}
where $Q_{-}:=B_o(\frac 12r)\times (s-\frac 34r^2,s-\frac 12r^2)$
and $Q_{+}:=B_o(\frac 12r)\times (s-\frac 14r^2,s)$, and where
$\bar{V}_f(Q_{-})$ denotes the volume of $Q_{-}$ with respect to
the space-time volume form $d\bar{\mu}$. Now we let $u$ be a solution,
which is not necessarily non-negative, and let $M_u$, $m_u$ be the
maximum and minimum of $u$ in $Q$. Similarly, let $M^{+}_u$, $m^{+}_u$ be the
maximum and minimum of $u$ in $Q_{+}$. Define
\[
\mu^{-}_u:=\frac{1}{\bar{\mu}(Q_{-})}\int_{Q_{-}}vd\bar{\mu},
\]
where $d\bar{\mu}$ denotes the natural product measure on $R\times M$:
$d\bar{\mu}=dt\times d\mu$, and where $d\mu=e^{-f}dv$.
Applying \eqref{Harineq} to the non-negative solutions $M_u-u$, $u-m_u$
yields
\[
M_u-\mu^{-}_u\leq e^{c(n,A)(1+Kr^2)}(M_u-M^{+}_u)
\]
and
\[
\mu^{-}_u-m_u\leq e^{c(n,A)(1+Kr^2)}(m^{+}_u-m_u),
\]
which imply that
\[
(M_u-m_u)\leq e^{c(n,A)(1+Kr^2)}(M_u-m_u)-e^{c(n,A)(1+Kr^2)}(M^{+}_u-m^{+}_u).
\]
If we define the oscillations
\[
\omega(u,Q):=M_u-m_u \quad \mathrm{and}\quad
\omega(u,Q_{+}):=M^{+}_u-m^{+}_u
\]
of $u$ over $Q$ and
$Q_{+}$, then
\begin{equation}\label{oscineq}
\omega(u,Q_{+})\leq\theta\omega(u,Q),
\end{equation}
where we assume $e^{c(n,A)(1+Kr^2)}>1$, and hence
$\theta=1-e^{-c(n,A)(1+Kr^2)}\in(0,1)$.

Now we consider $(y,t),(y',t')\in Q_\kappa$. Let
\[
\rho=2\max\{d(y,y'),\sqrt{t-t'}\}
\]
with $t\geq t'$. Then $(y',t')$ belongs to $Q_0:=B_y(\rho)\times(t-\rho^2,t)$.
We also define $\rho_i=2\rho_{i-1}$, $\rho_0=\rho$ and
$Q_i:=B_y(\rho_i)\times (t-\rho_i^2,t)$ for all $i\geq 1$. We easily
see that
\[
(Q_i)_+=Q_{i-1}.
\]
Hence, as long as $Q_i$ is contained in $Q$, \eqref{oscineq} yields
\[
\omega(u,Q_{i-1})\leq\theta\omega(u,Q_i)\quad\mathrm{and}\quad
\omega(u,Q_0)\leq\theta^i\omega(u,Q_i).
\]

Below, we consider two cases. If $\rho\leq (1-\kappa)r$, let $k$
be the integer such that
\[
2^k\leq (1-\kappa)r/\rho<2^{k+1}.
\]
Since $(y,t)\in Q_\kappa$, it follows that
\begin{equation*}
\begin{aligned}
Q_k&=B_y(2^k\rho)\times(t-4^k\rho^2,t)\\
&\subset B_y((1-\kappa)r)\times(t-(1-\kappa)^2r^2,t)\\
&\subset B_o(r)\times(s-r^2,s)=Q.
\end{aligned}
\end{equation*}
Hence we have
\[
\omega(u,Q_0)\leq \theta^k\omega(u,Q)
\leq \theta^{-1}(1-\kappa)^{-\alpha}\left(\frac{\rho}{r}\right)^\alpha\omega(u,Q)
\]
with $\alpha=-\log_2\theta$. This implies
\[
\frac{|u(y,t)-u(y',t')|}{[|t-t'|^{1/2}+d(y,y')]^\alpha}
\leq \frac{A_\kappa}{r^\alpha}\sup_Q\{|u|\}
\]
and conclusion follows, where $A_\kappa:=4\theta^{-1}(1-\kappa)^{-\alpha}$.
The second case is trivial. Indeed if $\rho>(1-\kappa)r$,
then the above inequality obviously holds. Therefore we
complete the proof of theorem.
\end{proof}

Using the Harnack inequality and the H\"{o}lder continuity
estimate, we immediately derive the following Liouville theorem.
\begin{corollary}\label{holliou}
Let $(M,g,e^{-f}dv)$ be an $n$-dimensional complete noncompact
smooth metric measure space. Assume that $Ric_f\geq 0$ and
$|f|(x)\leq A$ for some nonnegative constant $A$. Then any
solution $u$ of the $f$-harmonic equation which is bounded
from below (or above) is constant. Moreover, there exists
an $\alpha\in(0,1]$ such that any $f$-harmonic function
$u$ which satisfies
\begin{equation}\label{increcond}
\lim_{r\to\infty}\left(r^{-\alpha}\cdot\sup_{B_o(r)}\{|u|\}\right)=0
\end{equation}
for some fixed $o\in M$ is constant.
\end{corollary}
\begin{remark}
Corollary \ref{holliou} was also proved by Munteanu and Wang
\cite{[MuWa]}. We emphasize that this result can be regarded
as a direct consequence of Theorem \ref{THEO}. If $f$ is
constant and $\alpha=1$, then Corollary \ref{holliou} returns
to Cheng's Liouville property in \cite{[Cheng]}. If $f$ is
constant, this case appeared in \cite{[Saloff2]} (see also
Theorem 4.3 in \cite{[Saloff]}).
\end{remark}
\begin{proof}[Proof of Corollary \ref{holliou}]
We start to prove the first part of corollary. The conditions of
corollary imply the parabolic Harnack inequality (Corollary
\ref{coro1a}, $K=0$) and hence the corresponding elliptic
Harnack inequality. Assume that $u$ is a solution of the
$f$-harmonic equation which is bounded from below (if $u$ is
bounded from above, we then consider $-u$, which is still
bounded from below). Let
\[
m(u):=\inf_M\{u\}.
\]
Applying the elliptic Harnack inequality in the ball
$2B=B_o(2r)$ to the non-negative function
$v=u-m(u)$, we have that
\[
\sup_B\{u-m(u)\}\leq C(n,A)\cdot\inf_B\{u-m(u)\}.
\]
As the radius of $B=B_o(r)$ tend to infinity,
$\inf_B\{u-m(u)\}$ tends to zero. Therefore we conclude
that $u=m(u)$ is constant.

Below we will prove the second part of corollary. Because
$u$ has sublinear growth by condition \eqref{increcond},
then $\alpha$ can be taken in the interval $(0,1)$.
Let $\alpha$ be as given by Theorem \ref{MainHl}. Let
$u$ be a function satisfying $\Delta_fu=0$ and condition
\eqref{increcond}. Fix some $x\in M$ and $y$ such that
$d(x,y)\leq 1$. Applying Theorem \ref{MainHl} to $u$
in a ball $B_R=B_o(R)$ with $R$ so large that
$x,y\in \frac12B_R$, we find that
\begin{equation}\label{holdese}
|u(x)-u(y)|\leq \frac{C}{R^\alpha}\sup_{B_R}\{|u|\},
\end{equation}
where constant $C$ is independent of $R$. Since the above
inequality holds for all $R$ large enough, we can let $R$
tend to infinity to obtain that $|u(x)-u(y)|=0$. Since
$x,y\in M$ with $d(x,y)\leq 1$ are arbitrary and $M$ is
connected, we conclude that $u$ must be constant.
\end{proof}

\section{Two-sided Gaussian bounds on $f$-heat kernel}\label{sec2c}
In this section, we shall obtain upper and lower bound estimates
for the $f$-heat kernel on complete noncompact metric
measure space. The proof seems to be different from the classical
discussion of Li and Yau in \cite{[Li-Yau1]}. Our argument is
similar to the discussion in Grigor'yan \cite{[Grig]} and
Saloff-Coste \cite{[Saloff]}.

First, we show that the local $f$-Neumann Poincar\'e inequality
and the the local $f$-volume
doubling property imply a lower bound on the $f$-heat kernel. To
achieve this, we begin with by the following important lemma.
\begin{lemma}\label{lemm3.1}
Under the same assumptions of Theorem \ref{PoinDouHarn},
there exists a constant $c_4:=c_4(n,A)$ such that, for any
$x,y\in B_o(\frac 12R)$, and any $0<s<t<\infty$ and any
non-negative solution $u$ of the $f$-heat equation in
$M\times(0,\infty)$,
\[
\ln\left(\frac{u(x,s)}{u(y,t)}\right)\leq c_4
\left[\left(K+\frac{1}{R^2}+\frac{1}{s}\right)(t-s)+\frac{d^2(x,y)}{t-s}\right].
\]
\end{lemma}
\begin{proof}[Sketch proof of Lemma \ref{lemm3.1}]
Since Theorem \ref{PoinDouHarn} implies a parabolic Harnack inequality
of the $f$-heat equation, it is sufficient to prove the above
inequalities by carefully choosing different space-time solutions.
Please see Corollary 5.4.4 in \cite{[Saloff2]} or Corollary 5.4
in \cite{[Saloff1]} for a detailed proof.
\end{proof}

Using Lemma \ref{lemm3.1}, we can get a lower bound on the
$f$-heat kernel on complete metric measure spaces.
\begin{proposition}\label{prop3.1}
Under the same assumptions of Theorem \ref{PoinDouHarn},
there exists a constant $c_5:=c_5(n,A)$ such that, for any
$x,y\in B_o(\frac 12R)$ and any $0<t<\infty$, the $f$-heat
kernel $H(x,y,t)$ satisfies
\begin{equation}\label{lowheax}
H(x,y,t)\geq H(x,x,t)\exp
\left[-c_5\left(1+\frac{t}{R^2}+Kt+\frac{d^2(x,y)}{t}\right)\right].
\end{equation}
Moreover, there exists a constant $c_6:=c_6(n,A)$  such that, for any
$x,y\in B_o(\frac 12R)$ and any $0<t<R^2$
\begin{equation}\label{lowhe2x}
H(x,y,t)\geq \frac{e^{-c_6(1+Kt)}}{V_f(B_x(\sqrt{t}))}
\exp\left(-c_5\frac{d^2(x,y)}{t}\right).
\end{equation}
\end{proposition}

\begin{proof}
The proof follows from that of Theorem 5.4.11 in \cite{[Saloff2]}
with minor modifications. In fact using Lemma \ref{lemm3.1}, we
let $u(y,t)=H(x,y,t)$ with $x$ fixed and $s=t/2$ and then we get
\eqref{lowheax}, where we used the fact that $H(x,x,t)$ is
non-increasing.

Below we prove \eqref{lowhe2x}. Note that the conditions
of the proposition imply a parabolic Harnack inequality, which
leads to the on-diagonal $f$-heat kernel lower bound
\begin{equation}\label{hklb}
H(x,x,t)\geq e^{-c(n,A)(1+Kt)}\cdot V_f^{-1}(B_x(\sqrt{t}))
\end{equation}
for all $x\in M$ and $0<t<R^2$. Indeed we fix $0<t<R$ and
consider $\phi$ be a smooth function such that
$0\leq\phi\leq1$, $\phi=1$
on $B:=B_x(\sqrt{t})$ and $\phi=0$ on $M\setminus2B$. Define
\begin{equation*}
u(y,t)=\left\{ \begin{aligned}
         P_t\phi(y)\quad\mathrm{if}\quad t>0 \\
                  \phi(y)\quad \mathrm{if} \quad t\leq 0,
                          \end{aligned} \right.
                          \end{equation*}
$P_t=e^{t\Delta_f}$ be the heat semigroup of $\Delta_f$
on $L^2(M,\mu)$. Obviously, $u(y,t)$ satisfies
$(\partial_t-\Delta_f)u=0$ on $B\times(-\infty,\infty)$.
Applying the parabolic Harnack inequality, first to $u$, and
then to the $f$-heat kernel $(y,s)\to H(x,y,s)$, we have
\begin{equation*}
\begin{aligned}
1=u(x,0)&\leq e^{c(1+Kt)} u(x,t/2)\\
&=e^{c(1+Kt)}\int_{B(x,\sqrt{t})} H(x,y,t/2)\phi(y)d\mu(y)\\
&\leq e^{c(1+Kt)}\int_{B(x,2\sqrt{t})} H(x,y,t/2)d\mu(y)\\
&\leq e^{2c(1+Kt)}V_f(B_x(2\sqrt{t}))H(x,x,t)\\
&\leq e^{2c(1+Kt)}V_f(B_x(\sqrt{t}))2^{n+4A}e^{C(n,A)\sqrt{Kt}}H(x,x,t),
\end{aligned}
\end{equation*}
where in the last inequality we used \eqref{voldop}.
This gives \eqref{hklb} as desired. Hence \eqref{lowhe2x} then
easily follows by \eqref{lowheax} and \eqref{hklb}.
\end{proof}

Secondly, we can show that the local $f$-Neumann Poincar\'e inequality
and the local $f$-volume doubling property also imply an upper bound on
the $f$-heat kernel. To achieve this, the following integral
estimate is critically useful due to Davies \cite{[Davies]}.
\begin{lemma}[Davies \cite{[Davies]}]\label{lemm3.3}
Let $(M,g,e^{-f}dv)$ be an $n$-dimensional complete smooth
metric measure space. Let $\lambda_1>0$ be the bottom of the
$L^2$-spectrum of the $f$-Laplacian. Assume that $B_1$ and
$B_2$ are bounded subsets of $M$. Then
\[
\int_{B_1}\int_{B_2}H(x,y,t)d\mu(y)d\mu(x)\leq e^{-\lambda_1t}
V_f(B_1)^{1/2}V_f(B_2)^{1/2}\exp\left(-\frac{d^2(B_1,B_2)}{4t}\right),
\]
where $d(B_1,B_2)$ denotes the distance between the sets $B_1$ and $B_2$.
\end{lemma}

We now give an upper bound on the fundamental solution of
the $f$-heat equation.
\begin{proposition}\label{prop3.2}
Under the same assumptions of Theorem \ref{PoinDouHarn},
there exist constants $c_7$ and $c_8$ such that, for any
$x,y\in B_o(\frac 12R)$ and $0<t<R^2/4$, the $f$-heat
kernel $H(x,y,t)$ satisfies
\begin{equation}\label{lowhe2}
H(x,y,t)\leq \frac{e^{c_8(1+Kt)}}
{V_f(B_x(\sqrt{t}))}\exp\left(-c_7\frac{d^2(x,y)}{t}\right).
\end{equation}
\end{proposition}
\begin{proof}
Fix a fixed $y\in B_o(r)$ and $\delta>0$, applying Lemma \ref{lemm3.1}
to the positive solution $u(x,t)=H(x,y,t)$ by taking $s=t$ and
$t=(1+\delta)t$,
\[
H(x,y,t)\leq H(x',y,(1+\delta)t)\cdot\exp
\left\{c_4\left[\left(K+\frac{1}{R^2}+\frac{1}{t}\right)\delta t
+\frac{d^2(x,x')}{\delta t}\right]\right\}.
\]
Integrating over $x'\in B_x(\sqrt {t})$ gives
\begin{equation}\label{heaupp}
\begin{aligned}
H(x,y,t)&\leq \exp\left[c_4\left((K+R^{-2})\delta t
+\delta+\frac{1}{\delta}\right)\right]V^{-1}_f(B_x(\sqrt {t}))\\
&\quad\times\int_{B_x(\sqrt {t})}H(x',y,(1+\delta)t)d\mu(x').
\end{aligned}
\end{equation}
Applying Lemma \ref{lemm3.1} and the same argument to the
positive solution
\[
u(y,t)=\int_{B_x(\sqrt {t})}H(x',y,t)d\mu(x'),
\]
by taking $s=(1+\delta)t$ and $t=(1+2\delta)t$, we obtain
\begin{equation*}
\begin{aligned}
\int_{B_x(\sqrt {t})}H(x',y,(1+\delta)t)d\mu(x')\leq
\exp\left[c_4\left((K+R^{-2})\delta t+\delta+\frac{1}{\delta}\right)\right]
V^{-1}_f(B_y(\sqrt {t}))\\
\times\int_{B_y(\sqrt {t})}\int_{B_x(\sqrt {t})}H(x',y',(1+2\delta)t)d\mu(x')d\mu(y').
\end{aligned}
\end{equation*}
Substituting this into \eqref{heaupp} yields
\begin{equation*}
\begin{aligned}
H(x,y,t)&\leq \exp\left[2c_4\left((K+R^{-2})\delta t+\delta+\frac{1}{\delta}\right)\right]
V^{-1}_f(B_x(\sqrt {t}))V^{-1}_f(B_y(\sqrt {t}))\\
&\quad\times\int_{B_y(\sqrt {t})}
\int_{B_x(\sqrt {t})}H(x',y',(1+2\delta)t)d\mu(x')d\mu(y').
\end{aligned}
\end{equation*}
Combining this with Lemma \ref{lemm3.3}, we have
\begin{equation}\label{heaupp2}
\begin{aligned}
H(x,y,t)&\leq \exp\left[2c_4\left((K+R^{-2})\delta t+\delta
+\frac{1}{\delta}\right)-\lambda_1t\right]\\
&\quad\times V^{-1/2}_f(B_x(\sqrt {t}))V^{-1/2}_f(B_y(\sqrt {t}))
\exp\left[-\frac{d^2(B_x(\sqrt {t}),B_y(\sqrt {t}))}{4(1+2\delta)t}\right].
\end{aligned}
\end{equation}
Notice that if $d(x,y)\leq 2\sqrt{t}$, then $d(B_x(\sqrt {t}),B_y(\sqrt {t}))=0$ and
hence
\[
-\frac{d^2(B_x(\sqrt {t}),B_y(\sqrt {t}))}{4(1+2\delta)t}=0\leq 1-\frac{d^2(x,y)}{4(1+2\delta)t},
\]
and if $d(x,y)>2\sqrt{t}$, then $d(B_x(\sqrt {t}),B_y(\sqrt {t}))=d(x,y)-2\sqrt{t}$
hence
\[
-\frac{d^2(B_x(\sqrt {t}),B_y(\sqrt {t}))}{4(1+2\delta)t}
=-\frac{(d(x,y)-2\sqrt{t})^2}{4(1+2\delta)t}
\leq -\frac{d^2(x,y)}{4(1+2\delta)t}+\frac{1}{2\delta}.
\]
Therefore in any case, \eqref{heaupp2} becomes
\begin{equation}\label{heaupp3}
\begin{aligned}
H(x,y,t)&\leq\exp\left[1+2\left(c_4+\frac 14\right)\left((K+R^{-2})\delta t
+\delta+\frac{1}{\delta}\right)-\lambda_1t\right]\\
&\quad\times V^{-1/2}_f(B_x(\sqrt {t}))V^{-1/2}_f(B_y(\sqrt {t}))\exp\left(-\frac{d^2(x,y)}{4(1+2\delta)t}\right).
\end{aligned}
\end{equation}
Now we want to estimate
$(K+R^{-2})\delta t+\delta+\frac{1}{\delta}$ in \eqref{heaupp3}. Let
\[
\delta=\min\left\{\epsilon,\left[(K+R^{-2})t \right]^{-1/2}\right\}.
\]
If $\left[(K+R^{-2})t \right]^{-1/2}\leq \epsilon$, then
\[
(K+R^{-2})\delta t+\delta+\frac{1}{\delta}\leq2\left[(K+R^{-2})t \right]^{1/2}+\epsilon.
\]
If $\left[(K+R^{-2})t\right]^{-1/2}>\epsilon$, then we have
\begin{equation*}
\begin{aligned}
(K+R^{-2})\delta t+\delta+\frac{1}{\delta}&\leq\left[(K+R^{-2})t \right]
\epsilon+\epsilon+\frac{1}{\epsilon}\\
&\leq\left[(K+R^{-2})t \right]^{1/2}+\epsilon+\frac{1}{\epsilon}.
\end{aligned}
\end{equation*}
Hence, in either case, the right hand side of \eqref{heaupp3} can
be estimate by
\begin{equation}\label{heaupp4}
\begin{aligned}
H(x,y,t)&\leq \exp\left[1+2\left(c_4+\frac 14\right)
\left(2\left[(K+R^{-2})t \right]^{1/2}
+\epsilon+\frac{1}{\epsilon}\right)-\lambda_1t\right]\\
&\quad\times V^{-1/2}_f(B_x(\sqrt {t}))V^{-1/2}_f(B_y(\sqrt {t}))\exp\left(-\frac{d^2(x,y)}{4(1+2\epsilon)t}\right).
\end{aligned}
\end{equation}
Moreover the volume doubling property implies (see, e.g., Lemma 5.2.7
in \cite{[Saloff2]}) that
\begin{equation*}
\begin{aligned}
V_f(x,\sqrt{t})&\leq C(n,A)\exp\left(C(n,A)\sqrt{Kt}
\cdot\frac{d(x,y)}{\sqrt{t}}\right)V_f(y,\sqrt{t})\\
&\leq C(n,A)\exp\left(\bar{C}(n,A,\epsilon)Kt+\frac{d^2(x,y)}{8(1+2\epsilon)t}\right)
V_f(y,\sqrt{t}).
\end{aligned}
\end{equation*}
Substituting this into \eqref{heaupp4} and using $0<t<R^2/4$,
then the theorem follows.
\end{proof}

Combining Lemmas \ref{comp}, \ref{NeuPoin} and Propositions
\ref{prop3.1}, \ref{prop3.2} immediately yields two-sided
$f$-heat kernel bounds on complete noncompact metric
measure spaces.
\begin{theorem}\label{NeuSob2}
Let $(M,g,e^{-f}dv)$ be an $n$-dimensional complete noncompact
smooth metric measure space. If $Ric_f\geq-(n-1)K$ and
$|f|(x)\leq A$ on $B_o(2R)$ for some nonnegative constants
$K$ and $A$, then there exist positive constants $c_i$,
$i=5,6,7,8$, depending on $n$ and $A$ such that the $f$-heat
kernel $H(x,y,t)$ satisfies
\[
\frac{e^{-c_6(1+Kt)}}{V_f(x,\sqrt{t})}\exp\left(-c_5\frac{d^2(x,y)}{t}\right)
\leq H(x,y,t)\leq \frac{e^{c_8(1+Kt)}}{V_f(x,\sqrt{t})}
\exp\left(-c_7\frac{d^2(x,y)}{t}\right)
\]
for any $x,y\in B_o(R/2)$ and $0<t<R^2/4$.
\end{theorem}
\begin{remark}
In \cite{[Saloff3]} and \cite{[Grig3]}, Saloff-Coste and
Grigor'yan have proved that the conjunction of the
$f$-volume doubling property and the $f$-Neumann Poincar\'e
inequality is equivalent to the two-sided $f$-heat kernel bounds,
whereas we give concrete weighted curvature condition to achieve
these estimates.
\end{remark}

\section{$f$-Mean value inequality}\label{sec3}

In this section, the main objective is to derive a mean
value inequality on complete noncompact metric measure space,
which is a natural generalization of the Li-Schoen's result in
\cite{[Li-Sch]}. First, we give the following Poincar\'e inequality.
\begin{theorem}\label{mainh20}
Let $(M,g,e^{-f}dv)$ be a complete noncompact smooth metric measure
space. Let $o\in M$ and $R>0$. If $Ric_f\geq-(n-1)K$ and
$|f|(x)\leq A$ for some nonnegative constants $K$ and $A$,
then for any $\alpha\geq 1$, there exists constants $C_3$
and $C_4$ depending only on $\alpha$, $n$ and $A$ such that
\[
\int_{B_o(R)}|\phi|^\alpha d\mu\leq C_3\left(\frac{R}{1+\sqrt{K}R}\right)^\alpha
e^{C_4(1+\sqrt{K}R)}
\int_{B_o(R)}|\nabla \phi|^\alpha d\mu
\]
for any compactly supported function $\phi$ on $B_o(R)$.
In particular, the first Dirichlet eigenvalue $\mu_1$ of $f$-Laplacian
on $B_o(R)$ satisfies
\[
\mu_1\geq C^{-1}_3\left(\frac{R}{1+\sqrt{K}R}\right)^{-2}e^{-C_4(1+\sqrt{K}R)}.
\]
\end{theorem}
\begin{proof}[Sketch proof of Lemma \ref{mainh20}]
The proof is exactly the same as that of Corollary 1.1 proved by
Li-Schoen \cite{[Li-Sch]} except that the classical Laplacian
comparison is replaced by the generalized Laplacian comparison
(see Lemma \ref{comp})
\begin{equation*}
\begin{aligned}
\Delta_fr(x)&\leq(n-1+4A)\sqrt{K}\coth \sqrt{K}r\\
&\leq\frac{n-1+4A}{r}+(n-1+4A)\sqrt{K}.
\end{aligned}
\end{equation*}
Besides this, all the integration calculations should be done
with respect to the new measure $\mu$. To save the length of
paper, we omit details of the proof.
\end{proof}

We now proceed to derive the $L^2$ $f$-mean value inequality
by Theorem \ref{mainh20}, which is a weighted version of
Li-Schoen's result in \cite{[Li-Sch]}.
\begin{theorem}\label{mainh4}
Let $(M,g,e^{-f}dv)$ be a complete noncompact smooth metric
measure space. Assume that $Ric_f\geq-(n-1)K$ with $|f|(x)\leq A$
for some nonnegative constants $K$ and $A$. Let $o\in M$ and $R>0$,
and let $u$ be a
nonnegative $f$-subharmonic function defined on $B_o(R)$.
There exists a constant $C_5$, depending only on $n$ and $A$
such that for any $\tau\in(0, 1/2)$ we have
\[
\sup_{B_o((1-\tau)R)}u^2\leq\tau^{-C_5(1+\sqrt{K}R)}
V_f^{-1}(B_o(R))\int_{B_o(R)}u^2 d\mu.
\]
\end{theorem}
\begin{proof}
The proof is similar to the Li-Schoen's proof of Theorem 1.2
in \cite{[Li-Sch]}. We include it here for the reader's convenience.
Let $h$ be a harmonic function on $B_o((1-2^{-1}\tau)R)$ obtained
by the solving the Dirichlet boundary problem
\[
\Delta_fh=0\quad \mathrm{on}\quad B_o((1-\tau/2)R),
\]
and
\[
h=u \quad \mathrm{on} \quad \partial B_o((1-\tau/2)R).
\]
Since $u$ is nonnegative, by the maximum principle, the function
$h$ is positive on the ball $B_o((1-2^{-1}\tau)R)$. Moreover,
\[
u\leq h \quad \mathrm{on}\quad B_o((1-\tau/2)R).
\]
Using Lemmas \ref{comp}, \ref{NeuPoin} and \ref{NeuSob},
by the Moser iteration argument as in \cite{[MuWa]}, we have
the following elliptic Harnack inequality
\[
\sup_{B_o((1-\tau)R)}h\leq e^{c(n,A)(1+\sqrt{K}R)}\inf_{B_o((1-\tau)R)}h,
\]
where $c$ depends only on $n$ and $A$. In particular,
\begin{equation}
\begin{aligned}\label{guaan1}
\sup_{B_o((1-\tau)R)}u^2&\leq\sup_{B_o((1-\tau)R)}h^2\\
&\leq e^{c(n,A)(1+\sqrt{K}R)}\inf_{B_o((1-\tau)R)}h^2\\
&\leq e^{c(n,A)(1+\sqrt{K}R)}V_f^{-1}(B_o((1-\tau)R))\int_{B_o((1-\tau)R)}h^2 d\mu.
\end{aligned}
\end{equation}
Below we will estimate the $L^2(\mu)$-norm of $h$ in terms of the
$L^2(\mu)$-norm of $u$. By the triangle inequality, we have
\begin{equation}
\begin{aligned}\label{guaan2}
\int_{B_o((1-\tau)R)}h^2 d\mu
\leq&2\int_{B_o((1-\tau)R)}(h-u)^2 d\mu+2\int_{B_o((1-\tau)R)}u^2 d\mu\\
\leq&2\int_{B_o((1-\tau/2)R)}(h-u)^2 d\mu+2\int_{B_o(R)}u^2 d\mu.
\end{aligned}
\end{equation}
Since $(h-u)$ vanishes on $\partial B_o((1-\tau/2)R)$ we can apply
Theorem \ref{mainh20} to show that
\begin{equation*}
\begin{aligned}
\int_{B_o((1-\tau/2)R)}(h-u)^2 d\mu
\leq&\frac{C_3R^2}{\left(1+\sqrt{K}R\right)^2}e^{C_4(1+\sqrt{K}R)}
\int_{B_o((1-\tau/2)R)}|\nabla (h-u)|^2 d\mu\\
\leq&\frac{C_3R^2e^{C_4(1+\sqrt{K}R)}}{\left(1+\sqrt{K}R\right)^2}
\int_{B_o((1-\tau/2)R)}2(|\nabla h|^2+|\nabla u|^2)d\mu,
\end{aligned}
\end{equation*}
where we have used the triangle inequality again. Since the
Dirichlet integral of $h$ is least among all functions
which coincide with $h$ on the boundary, from above we conclude
that
\begin{equation}\label{guaan3}
\int_{B_o((1-\tau/2)R)}(h-u)^2 d\mu
\leq\frac{4C_3R^2e^{C_4(1+\sqrt{K}R)}}{\left(1+\sqrt{K}R\right)^2}
\int_{B_o((1-\tau/2)R)}|\nabla u|^2 d\mu.
\end{equation}
Now we use the fact that $u$ is $f$-subharmonic to estimate the
Dirichlet integral of $u$ in terms of the $L^2$-norm of $u$.
We have for any $\phi$ with compact support in $B_o(R)$
\begin{equation*}
\begin{aligned}
0&\leq\int_{B_o(R)}\phi^2u\Delta_fud\mu\\
&=-\int_{B_o(R)}\phi^2|\nabla u|^2d\mu
+2\int_{B_o(R)}\phi u\langle\nabla\phi,\nabla u\rangle d\mu\\
&\leq-\int_{B_o(R)}\phi^2|\nabla u|^2d\mu
+2\left(\int_{B_o(R)}\phi^2 |\nabla u|^2d\mu\right)^{1/2}
\left(\int_{B_o(R)}u^2|\nabla \phi|^2d\mu\right)^{1/2}.
\end{aligned}
\end{equation*}
Thus
\[
\int_{B_o(R)}\phi^2 |\nabla u|^2d\mu\leq 4\int_{B_o(R)}u^2|\nabla \phi|^2d\mu.
\]
We let $\phi(r(x))$ be a cut-off function given by a function of
$r(x)=r(o,x)$ alone, such that $\phi(r)=1$ on $B_o((1-\tau/2)R)$,
$\phi(r)=0$ on $\partial B_o(R)$, and satisfying
\[
|\nabla \phi|\leq \frac{c}{\tau R}.
\]
Then the above inequality becomes
\[
\int_{B_o((1-\tau/2)R)}|\nabla u|^2d\mu\leq\frac{4c^2}{\tau^2 R^2}
\int_{B_o(R)}u^2d\mu.
\]
Combining this with \eqref{guaan1}, \eqref{guaan2} and \eqref{guaan3}
yields
\begin{equation}
\begin{aligned}\label{guank}
\sup_{B_o((1-\tau)R)}u^2\leq&C
\left(\frac{32c^2C_3\tau^{-2}e^{C_4(1+\sqrt{K}R)}}
{(1+\sqrt{K}R)^2}+2\right)
V_f^{-1}(B_o((1-\tau)R))\int_{B_o(R)}u^2 d\mu\\
\leq&C_6\tau^{-C_7(1+\sqrt{K}R)}e^{C_8(1+\sqrt{K}R)}
V_f^{-1}(B_o((1-\tau)R))\int_{B_o(R)}u^2 d\mu
\end{aligned}
\end{equation}
for some new constants $C_i=C_i(n,A)$, $i=6,7,8$. To
finish the proof, we also need to estimate the $f$-volume of
$B_o(R)$ in terms of the volume of $B_o((1-\tau)R)$.
Recall the bound for $\Delta_fr^2$:
\[
\Delta_fr^2\leq 2(n+4A)+2\sqrt{K}(n-1+4A)r,
\]
and hence
\[
\int_{B_o(t)}\Delta_fr^2d\mu\leq 2(n+4A)V_f(t)
+2\sqrt{K}(n-1+4A)\int_{B_o(t)}rd\mu,
\]
where $V_f(t)=Vol_f(B_o(t))$. By Green formula, since
\[
\int_{B_o(t)}\Delta_fr^2d\mu
=\int_{\partial B_o(t)}\frac{\partial r^2}{\partial r}d\sigma
=2t\frac{\partial V_f(B_o(t))}{\partial t},
\]
then
\[
tV^{'}_f(t)\leq(n+4A)V_f(t)+\sqrt{K}(n-1+4A)tV_f(t).
\]
Hence the function $t^{-(n+4A)}e^{-\sqrt{K}(n-1+4A)t}V_f(t)$
is decreasing in $t\geq 0$. Therefore
\begin{equation}\label{rela}
V^{-1}_f(B_o((1-\tau)R))\leq V^{-1}_f(B_o(R))
\left(\frac{1}{1-\tau}\right)^{n+4A}\cdot e^{\sqrt{K} R\tau(n-1+4A)},
\end{equation}
where $0<\tau<1/2$. Combining this with \eqref{guank} completes
the proof of theorem.
\end{proof}

In the following, we show that the $L^p$ $f$-mean value inequality
for any $p\in(0,2]$ is a formal consequence of that given in
Theorem \ref{mainh4}.
\begin{theorem}\label{mainh5}
Under the same assumption of Theorem \ref{mainh4}, for any $p\in(0, 2]$,
there exists a constant $c$ depending only on $n$, $p$,
and $A$ such that
\[
\sup_{B_o((1-\tau)R)}u^p\leq\tau^{-c(1+\sqrt{K}R)}
V^{-1}_f(R)\int_{B_o(R)}u^pd\mu
\]
for any $\tau\in(0, 1/2)$, where $V^{-1}_f(R):=V^{-1}_f(B_o(R))$.
\end{theorem}
\begin{proof}
The proof is similar to the proof of Theorem 2.1 in \cite{[Li-Sch]}.
However, for the sake of completeness, we include the details
here. By Theorem \ref{mainh4}, for any $\delta\in(0, 1/2]$,
$\theta\in[1/2,1-\delta]$, we have
\[
\sup_{B_o(\theta R)}u^2\leq\delta^{-C_5(1+\sqrt{K}R)}
V_f^{-1}((\theta+\delta)R)\int_{B_o((\theta+\delta)R)}u^2 d\mu.
\]
Since $\theta+\delta\geq1/2$, this inequality implies
\[
\sup_{B_o(\theta R)}u^2\leq\delta^{-C_5(1+\sqrt{K}R)}
V_f^{-1}(2^{-1}R)\int_{B_o((\theta+\delta)R)}u^2 d\mu.
\]
We also note that
\[
\int_{B_o((\theta+\delta)R)}u^2 d\mu\leq
\left(\sup_{B_o((\theta+\delta)R)}u^2\right)^{1-p/2}
\int_{B_o((\theta+\delta)R)}u^p d\mu.
\]
Hence
\[
\sup_{B_o(\theta R)}u^2\leq\delta^{-C_5(1+\sqrt{K}R)}
V_f^{-1}(2^{-1}R)\left(\sup_{B_o((\theta+\delta)R)}u^2\right)^{1-p/2}
\int_{B_o((\theta+\delta)R)}u^p d\mu.
\]
If we set
\[
M(\theta):=\sup_{B_o(\theta R)}u^2
\]
and
\[
N:=V_f^{-1}(2^{-1}R)\int_{B_o(R)}u^pd\mu,
\]
we have shown
\[
M(\theta)\leq N\delta^{-C_5(1+\sqrt{K}R)}M(\theta+\delta)^{1-p/2}
\]
for any $\delta\in(0, 1/2]$ and $\theta\in[1/2,1-\delta]$.
Choosing
\[
\theta_0=1-\tau \quad\mathrm{and}\quad\theta_i=\theta_{i-1}+2^{-i}\tau
\]
for $i=1,2,3,...$, we have that
\[
M(\theta_{i-1})\leq N_12^{iC_5(1+\sqrt{K}R)}M(\theta_i)^\lambda,
\]
where $\lambda=1-p/2$ and $N_1=N\tau^{-C_5(1+\sqrt{K}R)}$. Iterating yields
\[
M(\theta_0)\leq K_1^{\Sigma^j_{i=1}\lambda^{i-1}}
2^{C_5(1+\sqrt{K}R)\Sigma^j_{i=1}i\lambda^{i-1}}M(\theta_j)^{\lambda^j}
\]
for any $j\geq1$. Letting $j$ tend to infinity yields
\[
M(\theta_0)\leq\tau^{-C_9(1+\sqrt{K}R)}
\left[V_f^{-1}(2^{-1}R)\int_{B_o(R)}u^pd\mu\right]^{2/p},
\]
where $C_9$ depends only on $n$, $p$ and $A$. By the definition of
$M(\theta_0)$, we have
\[
\sup_{B_o((1-\tau)R)}u^p\leq
\tau^{-2^{-1}pC_9(1+\sqrt{K}R)}V_f^{-1}(2^{-1}R)\int_{B_o(R)}u^pd\mu.
\]
Finally, by the relation \eqref{rela}, i.e.,
\[
V^{-1}_f(2^{-1}R)\leq C(n,A)e^{C(1+\sqrt{K}R)}V^{-1}_f(R),
\]
the theorem follows.
\end{proof}

\section{$L^p$-Liouville theorem}\label{sec4}
In this section, we will study various $L^p$-Liouville theorems
on complete noncompact smooth metric measure spaces. Our results
extend the classical $L^p$-Liouville theorems obtained by Li
and Schoen in \cite{[Li-Sch]} and P. Li \cite{[Li0]} and their
weighted versions proved by X.-D. Li in \cite{[LD]}.
\subsection{The $0<p<1$ case}\label{sub6.2}
For $0<p<1$, we have a new weighted version of
Li-Schoen's $L^p$-Liouville theorem in \cite{[Li-Sch]}.
\begin{theorem}\label{liouv2}
Let $(M,g,e^{-f}dv)$ be an $n$-dimensional complete noncompact
smooth metric measure space. Assume that $f$ is bounded, and
there exists a constant $\delta(n)>0$ depending only on $n$,
such that, for some point $o\in M$, the
Bakry-\'{E}mery Ricci curvature satisfies
\[
Ric_f\geq-\delta(n)r^{-2}(x),
\]
whenever the distance from $o$ to $x$, $r(x)$, is sufficiently
large. Then any nonnegative $L^p(\mu)$-integrable ($0<p<1$)
$f$-subharmonic function must be identically zero.
\end{theorem}

\begin{proof}[Proof of Theorem \ref{liouv2}]
The proof is similar to the arguments of
Li and Schoen (see Theorem 2.5 in \cite{[Li-Sch]}).
Since the arguments leading to Theorem \ref{mainh5} are local,
by choosing more or less $\tau=1/2$, we have the following
$L^p$ $f$-mean value inequality
\begin{equation}\label{kbds}
\sup_{B_o(R/2)}u^p\leq \cdot2^{c(1+\sqrt{K(x,5R)}R)}
V_f^{-1}(B_o(R))\int_{B_o(R)}u^p d\mu
\end{equation}
for nonnegative $f$-subharmonic functions $u$ on $B_x(5R)$,
where $Ric_f\geq-(n-1)K(x,5R)$ and $|f|(x)\leq A$ for some
nonnegative constants $K$ and $A$ on $B_x(5R)$. Here the
constant $c$ depends on $n$, $p$ and $A$. In the following,
we will use \eqref{kbds} to show that $u$ must vanish at
infinity if the nonnegative function $u$ is $f$-subharmonic
on $M$ with $u\in L^p(\mu)$ ($0<p<1$). In fact, by the volume
comparison theorem mentioned above, under the hypothesis
on $Ric_f$ and $f$, $M$ must be of $f$-infinite volume
and $u$ must be identically zero.

Let $x\in M$ and consider a minimal geodesic $\gamma$ joining
$o$ to $x$ such that $\gamma(0)=o$ and $\gamma(T)=x$, where
$T=r(o,x)$. We then define a set of values $\{t_i\in[O,T]\}^k_{i=0}$
satisfying
\[
t_0=0,\quad t_1=1+\beta,\quad
\ldots,\quad t_i=2\sum^i_{j=0}\beta^j-1-\beta^i,
\]
where $\beta>1$ to be chosen later, and
$t_k=2\sum^k_{j=0}\beta^j-1-\beta^k$ is the largest such value
with $t_k<T$. We denote the points $x_i=\gamma(t_i)$ and they
obviously satisfy
\[
r(x_i,x_{i+1})=\beta^i+\beta^{i+1},\quad r(o,x_i)=t_i\quad
\mathrm{and}\quad r(x_k,x)<\beta^k+\beta^{k+1}.
\]
Moreover, the set of geodesic balls $B_{x_i}(\beta^i)$ cover
$\gamma([0,2\sum^k_{j=0}\beta^j-1])$ and they have disjoint
interiors. We now claim that
\begin{equation}\label{kbdenshi}
V_f(B_{x_k}(\beta^k))\geq C
\left(\frac{\beta^{n+4A}}{(\beta+2)^{n+4A}-\beta^{n+4A}}\right)^k
V_f(B_o(1))
\end{equation}
for a fixed $\beta>2/(2^{1/n}-1)^{-1}>1$.
The proof of this claim essentially follows the arguments of
Cheeger-Gromov-Taylor in \cite{[CGT]}. For the sake of completeness,
we will outline the proof of this claim again.

For each $1\leq i\leq k$, a relative comparison
theorem (see (4.10) in \cite{[WW]}) argument shows that
\begin{equation*}
\begin{aligned}
V_f(B_{x_i}(\beta^i))&\geq D_i\left[V_f(B_{x_i}(\beta^i+2\beta^{i-1}))
-V_f(B_{x_i}(\beta^i))\right]\\
&\geq D_iV_f(B_{x_{i-1}}(\beta^{i-1})),
\end{aligned}
\end{equation*}
where
\[
D_i=\frac{\int^{\beta^i\sqrt{K(x_i,\beta^i+2\beta^{i-1})}}_0\sinh^{n-1+4A} tdt}
{\int^{(\beta^i+2\beta^{i-1})\sqrt{K(x_i,\beta^i+2\beta^{i-1})}}_{
\beta^i\sqrt{K(x_i,\beta^i+2\beta^{i-1})}}\sinh ^{n-1+4A}tdt},
\]
since $Ric_f\geq-(n-1)K(x_i,\beta^i+2\beta^{i-1})$ and
$|f|(x)\leq A$ for some nonnegative constants $K$ and $A$ on
$B_{x_i}(\beta^i+2\beta^{i-1})$. Iterating this inequality,
we conclude that
\begin{equation}\label{sum}
V_f(B_{x_k}(\beta^k))\geq V_f(B_o(1))\prod^k_{i=1}D_i.
\end{equation}
Since $r(o,x_i)=2\sum^i_{j=0}\beta^j-1-\beta^i$, the curvature
assumption implies that
\begin{equation*}
\begin{aligned}
\sqrt{K(x_i,\beta^i+2\beta^{i-1})}
&\leq\sqrt{\delta(n)}\cdot\left(2\sum^{i-2}_{j=0}\beta^j-1\right)^{-1}\\
&=\sqrt{\delta(n)}\cdot
\frac{\beta-1}{2\beta^{i-1}-\beta-1}
\end{aligned}
\end{equation*}
for sufficiently large $i$. Hence
\begin{equation*}
\begin{aligned}
\beta^i\sqrt{K(x_i,\beta^i+2\beta^{i-1})}
&\leq\sqrt{\delta(n)}\cdot
\frac{(\beta-1)\beta^i}{2\beta^{i-1}-\beta-1}\\
&=\sqrt{\delta(n)}\cdot
\frac{(\beta-1)\beta}{2-\beta^{2-i}-\beta^{1-i}}
\end{aligned}
\end{equation*}
which can be made arbitrarily small for a fixed
$\beta>2/(2^{1/n}-1)^{-1}>1$ by choosing
$\delta(n)$ to be sufficiently small.
Hence $D_i$ has the
following approximation
\begin{equation*}
\begin{aligned}
D_i&\sim\frac{(\beta^i)^{n+4A}}{(\beta^i+2\beta^{i-1})^{n+4A}-(\beta^i)^{n+4A}}\\
&=\frac{\beta^{n+4A}}{(\beta+2)^{n+4A}-\beta^{n+4A}}
\end{aligned}
\end{equation*}
by simply approximating $\sinh t$ with $t$.
Hence \eqref{kbdenshi} follows by combining \eqref{sum}.

\vspace{0.5em}

In the following, we shall estimate $V_f(B_x(\beta^{k+1}))$. We achieve it by two
cases.

Case 1: $r(x,x_k)\leq \beta^k(\beta-1)$. In this case, we see that
\[
B_{x_k}(\beta^k)\subset B_x(\beta^{k+1}),
\]
and hence
\[
V_f(B_{x_k}(\beta^k))\leq V_f(B_x(\beta^{k+1})).
\]
Combining this with \eqref{kbdenshi}, we conclude that
\[
V_f(B_x(\beta^{k+1}))\geq C\left(\frac{\beta^{n+4A}}
{(\beta+2)^{n+4A}-\beta^{n+4A}}\right)^kV_f(B_o(1)).
\]

Case 2: $r(x,x_k)>\beta^k(\beta-1)$. In this setting, we see that
\[
B_{x_k}(\beta^k)\subset B_x\big(r(x,x_k)+\beta^k\big)
\backslash B_x\big(r(x,x_k)-\beta^k\big).
\]
Using a relative comparison theorem, we have that
\begin{equation*}
\begin{aligned}
V_f(B_x(\beta^k))&\geq D\left[V_f\big(B_x(r(x,x_k)+\beta^k)\big)
-V_f\big(B_x(r(x,x_k)-\beta^k)\big)\right]\\
&\geq D\cdot V_f(B_{x_k}(\beta^k)),
\end{aligned}
\end{equation*}
where
\[
D=\frac{\int^{\beta^k\sqrt{K(x,r(x,x_k)+\beta^k)}}_0\sinh^{n-1+4A} tdt}
{\int^{(r(x,x_k)+\beta^k)\sqrt{K(x,r(x,x_k)+\beta^k)}}_{
(r(x,x_k)-\beta^k)\sqrt{K(x,r(x,x_k)+\beta^k)}}\sinh^{n-1+4A} tdt}
\]
Argument as above, since
\begin{equation*}
\begin{aligned}
(r(x,x_k)+\beta^k)\sqrt{K(x,r(x,x_k)+\beta^k)}&\leq
(\beta^{k+1}+2\beta^k)\sqrt{K(x,r(x,x_k)+\beta^k)}\\
&\leq\frac{\sqrt{\delta(n)}}{2}\cdot\beta(\beta-1)
\end{aligned}
\end{equation*}
can be made sufficiently small, we can approximate $D$ by
\[
D\sim\frac{\beta^{n+4A}}{(\beta+2)^{n+4A}}.
\]
Combining this with \eqref{kbdenshi} yields
\begin{equation}
\begin{aligned}\label{guguanj}
V_f(B_x(\beta^{k+1}))&\geq\frac{C\beta^{n+4A}}{(\beta+2)^{n+4A}}
\left(\frac{\beta^{n+4A}}
{(\beta+2)^{n+4A}-\beta^{n+4A}}\right)^{k+1}V_f(B_o(1))\\
&\geq \tilde{C}\left(\frac{\beta^{n+4A}}
{(\beta+2)^{n+4A}-\beta^{n+4A}}\right)^kV_f(B_o(1)),
\end{aligned}
\end{equation}
where $\tilde{C}$ depends on $n$, $A$ and $\beta$. In any case,
\eqref{guguanj} is valid.

If we let $x\to \infty$, the value $k\to \infty$. Note that the choice of
$\beta$ ensures that
\[
\frac{\beta^{n+4A}}{(\beta+2)^{n+4A}-\beta^{n+4A}}>1,
\]
and hence the right hand side of \eqref{guguanj} tends to infinity.

On the other hand, let us now apply to any point $x$ sufficiently
far from $o$. The assumption of theorem asserts that
$R\sqrt{K(x,5R)}$ is bounded from above. Combining this fact with
\eqref{kbds}, we have
\begin{equation}\label{contr}
u^p(x)\leq CV^{-1}_f(B_o(R)),
\end{equation}
where $C$ also depends on the $L^p$-norm of $u$. Using
the value $R=\beta^{k+1}$ in \eqref{contr}, the right hand
side of \eqref{contr} vanishes as $x\to\infty$, thus proving
that $u(x)\to 0$ as $x\to\infty$ and Theorem \ref{liouv2}
follows by the maximum principle.
\end{proof}

\subsection{The $p=1$ case}\label{sub6.3}
The proof of this case is more complex. We shall follow
the arguments of Li \cite{[Li0]} (see also \cite{[Li2]}
and \cite{[LD]}) to derive the following result.
\begin{theorem}\label{liouL0}
Let $(M,g,e^{-f}dv)$ be an $n$-dimensional complete noncompact
smooth metric measure space. Assume that $f$ is bounded, and
there exists a constant $C>0$, such that, for some point $o\in M$,
the Bakry-\'{E}mery Ricci curvature satisfies
\[
Ric_f\geq-C(1+r^2(x)),
\]
where $r(x)$ denotes the distance from $o$ to $x$. Then any
nonnegative $L^1(\mu)$-integrable $f$-subharmonic function
must be identically constant.
\end{theorem}

Following the trick of P. Li \cite{[Li0]} (see also \cite{[LD]})
to prove Theorem \ref{liouL0}, at first, we need the following
integration by parts formula.
\begin{theorem}\label{liouL1}
Under the same assumptions of Theorem \ref{liouL0}, for any
nonnegative $L^1(\mu)$-integrable $f$-subharmonic function $g$,
we have
\begin{equation*}
\int_M {\Delta_f}_y H(x,y,t)g(y)d\mu(y)
=\int_M H(x,y,t)\Delta_fg(y)d\mu(y).
\end{equation*}
\end{theorem}

\begin{proof}[Proof of Theorem \ref{liouL1}]
Similar to the proof of Theorem 1 in \cite{[Li0]}
(see also Theorem 6.1 in \cite{[LD]}),
applying the Green formula on $B_o(R)$, we have
\begin{equation*}
\begin{aligned}
&\left|\int_{B_o(R)}{\Delta_f}_y H(x,y,t)g(y)d\mu(y)
-\int_{B_o(R)}H(x,y,t)\Delta_f g(y)d\mu(y)\right|\\
=&\left|\int_{\partial
B_o(R)}\frac{\partial}{\partial r} H(x,y,t)g(y)d\mu_{\sigma,R}(y)-\int_{\partial B_o(R)}H(x,y,t)\frac{\partial}{\partial r}g(y)d\mu_{\sigma,R}(y)\right|\\
\leq&\int_{\partial B_o(R)}|\nabla H|(x,y,t)g(y)d\mu_{\sigma,R}(y)
+\int_{\partial B_o(R)}H(x,y,t)|\nabla g|(y)d\mu_{\sigma,R}(y),
\end{aligned}
\end{equation*}
where $\mu_{\sigma,R}$ denotes the weighted area measure induced by $\mu$
on $\partial B_o(R)$. In the following we shall prove that the above two
boundary integrals vanish as $R\to \infty$, which can be achieved
by five steps.

\emph{Step 1}.
In Theorem \ref{mainh5}, we have show that
any nonnegative subharmonic function $g(x)$ must satisfy
\[
\sup_{B_o(R)}g(x)\leq e^{c(1+R\sqrt{K(R)})}V_f^{-1}(2R)\int_{B_o(2R)}g(y)d\mu(y)
\]
for some constant $c=c(n,A)$, where $-(n-1)K(R)$ is the
lower bound of the Bakry-\'{E}mery Ricci curvature on
$B_o(4R)$ and  $|f|\leq A$. Applying our
theorem assumption, we have the estimate
\begin{equation} \label{mjbd}
\sup_{B_o(R)}g(x)\leq Ce^{\alpha R^2}V_f^{-1}(2R)\|g\|_{L^1(\mu)}
\end{equation}
for some constants $\alpha:=\alpha(n,A)$ and $C:=C(n,A)$.
Consider $\phi(y)=\phi(r(y))$ to be a nonnegative cut-off function
such that $0\leq\phi\leq1$, $|\nabla\phi|\leq\sqrt{3}$ and
\begin{equation*}
\phi(r(y))=\left\{ \begin{aligned}
&1&&\mathrm{on}\,\,B_o(R+1)\backslash B_o(R),\\
&0&&\mathrm{on}\,\,B_o(R-1)\cup(M\backslash B_o(R+2)).\\
\end{aligned} \right.
\end{equation*}
Since $g$ is $f$-subharmonic function, by the Schwarz inequality we have
\begin{equation*}
\begin{aligned}
0\leq\int_M\phi^2g\Delta_fg d\mu
=&-\int_M \nabla(\phi^2g)\nabla g d\mu\\
=&-2\int_M \phi g\langle\nabla\phi\nabla g\rangle d\mu
-\int_M \phi^2|\nabla g|^2d\mu\\
\leq&2\int_M |\nabla\phi|^2g^2d\mu
-\frac12\int_M \phi^2|\nabla g|^2d\mu.
\end{aligned}
\end{equation*}
Then using the definition of $\phi$ and \eqref{mjbd}, we have that
\begin{equation*}
\begin{aligned}
\int_{B_o(R+1)\backslash B_o(R)}|\nabla g|^2d\mu\leq&
4\int_M|\nabla \phi|^2g^2d\mu
\leq12\int_{B_o(R+2)}g^2d\mu\\
\leq&12\sup_{B_o(R+2)}g\cdot\|g\|_{L^1(\mu)}\\
\leq&\frac{Ce^{\alpha(R+2)^2}}{V_f(2R+4)}
\cdot\|g\|_{L^1(\mu)}^2.
\end{aligned}
\end{equation*}
On the other hand, using the Schwarz inequality, we get
\begin{equation*}
\begin{aligned}
\int_{B_o(R+1)\backslash B_o(R)}|\nabla g|d\mu\leq&
\left(\int_{B_o(R+1)\backslash B_o(R)}|\nabla g|^2d\mu\right)^{1/2}
\cdot [V_f(R+1)\backslash V_f(R)]^{1/2}\\
\leq& \left(\int_{B_o(R+1)\backslash B_o(R)}|\nabla g|^2d\mu\right)^{1/2}
\cdot V_f(2R+4)^{1/2}.
\end{aligned}
\end{equation*}
Combining the above two inequalities, we have
\begin{equation}\label{guji}
\int_{B_o(R+1)\backslash B_o(R)}|\nabla g|d\mu\leq C_{10}e^{\alpha R^2}
\cdot\|g\|_{L^1(\mu)},
\end{equation}
where $C_{10}=C_{10}(n,A)$.\\

\emph{Step 2}. We first estimate the $f$-heat kernel $H(x,y,t)$.
Recall that, by Theorem \ref{NeuSob2}, the $f$-heat kernel $H(x,y,t)$ satisfies
\[
H(x,y,t)\leq \frac{e^{c_8(1+K(R)t)}}
{V_f(B_x(\sqrt{t}))}\exp\left(-c_7\frac{d^2(x,y)}{t}\right)
\]
for all $x,y\in B_o(R)$ and $0<t<R^2/8$, where $-(n-1)K(R)$ is the
lower bound of the Bakry-\'{E}mery Ricci curvature on $B_o(2R)$.
Here the constants $c_7$ and $c_8$ depending on $n$ and $A$.
Combining this with the assumption of our theorem, we deduce that
\begin{equation}\label{guji2}
H(x,y,t)\leq \frac{C}
{V_f(B_x(\sqrt{t}))}\exp\left(-c_7\frac{d^2(x,y)}{t}+\alpha R^2t\right)
\end{equation}
for all $x,y\in B_o(R)$ and $0<t<R^2/8$.
Then combining \eqref{guji} with \eqref{guji2} gives
\begin{equation*}
\begin{aligned}
J_1:=&\int_{B_o(R+1)\backslash B_o(R)}H(x,y,t)|\nabla g|(y)d\mu(y)\\
\leq&\left(\sup_{y\in {B_o(R+1)\backslash B_o(R)} }H(x,y,t)\right)
\int_{B_o(R+1)\backslash B_o(R)}|\nabla g|d\mu\\
\leq&\frac{C_{11}\|g\|_{L^1(\mu)}}{V_f(B_x(\sqrt{t}))}
\cdot\exp\left[-c_7\frac{(R-d(o,x))^2}{t}+\alpha R^2t+\alpha R^2\right],
\end{aligned}
\end{equation*}
where $C_{11}=C_{11}(n,A)$. Note that
\begin{equation*}
\begin{aligned}
-c_7\frac{(R-d(o,x))^2}{t}&+\alpha R^2t+\alpha R^2\\
&=\left(\alpha t+\alpha-\frac{c_7}{t}\right)R^2+\frac{2c_7}{t}Rd(o,x)-c_7\frac{d^2(o,x)}{t}\\
&\leq\left(\alpha t+\alpha-\frac{c_7}{t}\right)R^2+\frac{c_7}{2t}R^2+c_7\frac{d^2(o,x)}{t}\\
&=\left(\alpha t+\alpha-\frac{c_7}{2t}\right)R^2+c_7\frac{d^2(o,x)}{t}.
\end{aligned}
\end{equation*}
Thus, for $T$ sufficiently small and for all $t\in (0,T)$ there exists some fixed
constant $\beta>0$ such that
\[
J_1
\leq C_{11}\|g\|_{L^1(\mu)}V^{-1}_f(B_x(\sqrt{t}))
\cdot\exp\left(-\beta R^2t+c_7t^{-1}d^2(o,x)\right).
\]
Hence, for all $t\in(0,T)$ and all $x\in M$, $J_1$ tends to zero as $R$
tends to infinity.\\

\emph{Step 3}.
Below we shall estimate the gradient of $H$. Here we adapt the Li's proof
trick (see Section 18 in \cite{[Li2]}). Consider the integral
with respect to $d\mu$:
\begin{equation*}
\begin{aligned}
\int_M\phi^2(y)|\nabla H|^2(x,y,t)
&=-2\int_M\big\langle H(x,y,t)\nabla\phi(y),\phi(y)\nabla H(x,y,t)\big\rangle\\
&\quad-\int_M\phi^2(y)H(x,y,t)\Delta_f H(x,y,t)\\
&\leq2\int_M|\nabla\phi|^2(y)H^2(x,y,t)+\frac 12\int_M\phi^2(y)|\nabla H|^2(x,y,t)\\
&\quad-\int_M\phi^2(y)H(x,y,t)\Delta_f H(x,y,t).
\end{aligned}
\end{equation*}
This implies
\begin{equation}
\begin{aligned}\label{cutfuest}
&\int_{B_o(R+1)\backslash B_o(R)}|\nabla H|^2
\leq\int_M\phi^2(y)|\nabla H|^2(x,y,t)\\
&\leq4\int_M|\nabla\phi|^2H^2-2\int_M\phi^2H\Delta_f H\\
&\leq12\int_{B_o(R+2)\backslash B_o(R-1)}H^2+2\int_{B_o(R+2)\backslash B_o(R-1)}H|\Delta_f H|\\
&\leq12\int_{B_o(R+2)\backslash B_o(R-1)}H^2
+2\left(\int_{B_o(R+2)\backslash B_o(R-1)}H^2\right)^{\frac 12}
\left(\int_M(\Delta_f H)^2\right)^{\frac 12}.
\end{aligned}
\end{equation}
It is known that the heat semi-group is contractive in $L^1(\mu)$, hence
\[
\int_M H(x,y,t)d\mu(y)\leq 1.
\]
Using this and \eqref{guji2}, we can estimate
\begin{equation}
\begin{aligned}\label{cutfuest2}
\int_{B_o(R+2)\backslash B_o(R-1)}&H^2(x,y,t)d\mu
\leq \sup_{y\in B_o(R+2)\backslash B_o(R-1)}H(x,y,t)\\
&\leq \frac{C_{12}}{V_f(B_x(\sqrt{t}))}
\times\exp\left[-c_7\frac{(R-1-d(o,x))^2}{t}+\alpha R^2t\right].
\end{aligned}
\end{equation}
Also, we \emph{claim} that there exists a constant $C_{13}>0$ such that
\begin{equation}\label{Lapest}
\int_M(\Delta_f H)^2(x,y,t)d\mu\leq\frac{C_{13}}{t^2}H(x,x,t).
\end{equation}
To prove this inequality, we first derive the inequality for any
Dirichlet $f$-heat kernel $H$ defined on a compact subdomain of
$M$. Using the fact that $f$-heat kernel on $M$ can be obtained
by taking limits of Dirichlet $f$-heat kernels on a compact
exhaustion of $M$, then \eqref{Lapest} follows. Indeed, if
$H(x,y,t)$ is a Dirichlet $f$-heat kernel on a compact
subdomain $\Omega\subset M$, using the eigenfunction
expansion, then $H(x,y,t)$ can be written as the form
\[
H(x,y,t)=\sum^{\infty}_ie^{-\lambda_it}\psi_i(x)\psi_i(y),
\]
where $\{\psi_i\}$ are orthonormal basis of the space of $L^2(\mu)$ functions
with Dirichlet boundary value satisfying the equation
\[
\Delta_f\psi_i=-\lambda_i\psi_i.
\]
Differentiating with respect to the variable $y$, we have
\[
\Delta_fH(x,y,t)=-\sum^{\infty}_i\lambda_ie^{-\lambda_it}\psi_i(x)\psi_i(y).
\]
Noticing that $s^2e^{-2s}\leq C_{13}e^{-s}$ for all
$0\leq s<\infty$, therefore
\[
\int_M(\Delta_fH)^2d\mu(y)\leq C_{13}t^{-2}
\sum^{\infty}_ie^{-\lambda_it}\psi^2_i(x)
=C_{13}t^{-2} H(x,x,t)
\]
and claim \eqref{Lapest} follows. Now combining  \eqref{cutfuest},
\eqref{cutfuest2} and \eqref{Lapest}, we obtain
\begin{equation*}
\begin{aligned}
\int_{B_o(R+1)\backslash B_o(R)}|\nabla H|^2d\mu\leq&
C_{14}\left[V^{-1}_f+t^{-1}V^{-\frac 12}_f H^{\frac 12}(x,x,t)\right]\\ &\times\exp\left[-c_7\frac{(R-1-d(o,x))^2}{2t}+\alpha R^2t\right],
\end{aligned}
\end{equation*}
where $V_f:=V_f(B_x(\sqrt{t}))$. Applying Schwarz inequality, we have
\begin{equation}
\begin{aligned}\label{gujiest2}
\int_{B_o(R+1)\backslash B_o(R)}|\nabla H|d\mu&\leq
\left[V_f(B_o(R+1))\backslash V_f(B_o(R))\right]^{1/2}\\
&\quad\times\left[\int_{B_o(R+1)\backslash B_o(R)}|\nabla H|^2d\mu\right]^{1/2}\\
&\leq V^{1/2}_f(B_o(R+1))
\left[V^{-1}_f+t^{-1}V^{-\frac 12}_f H^{\frac 12}(x,x,t)\right]^{1/2}\\
&\quad\times\exp\left[-c_7\frac{(R-1-d(o,x))^2}{4t}+\frac{\alpha}{2} R^2t\right].
\end{aligned}
\end{equation}
Therefore, by \eqref{mjbd}, \eqref{gujiest2} and Schwarz inequality we see that
\begin{equation*}
\begin{aligned}
J_2:=&\int_{B_o(R+1)\backslash B_o(R)}|\nabla H(x,y,t)|g(y)d\mu(y)\\
\leq&\sup_{y\in B_o(R+1)\backslash B_o(R)}g(y)\cdot
\int_{B_o(R+1)\backslash B_o(R)}|\nabla H(x,y,t)|d\mu(y)\\
\leq& \frac{Ce^{\alpha (R+1)^2}\|g\|_{L^1(\mu)}}{V_f(B_o(2R+2))}
\cdot V^{1/2}_f(B_o(R+1))
\left[V^{-1}_f+t^{-1}V^{-\frac 12}_f H^{\frac 12}(x,x,t)\right]^{1/2}\\
&\quad\times\exp\left[-c_7\frac{(R-1-d(o,x))^2}{4t}+\frac{\alpha}{2} R^2t\right]\\
\leq& \frac{C_{15}\|g\|_{L^1(\mu)}}{V^{1/2}_f(B_o(2R+2))}
\cdot\left[V^{-1}_f+t^{-1}V^{-\frac 12}_f H^{\frac 12}(x,x,t)\right]^{1/2}\\
&\quad\times\exp\left[-c_7\frac{(R-1-d(o,x))^2}{4t}+\frac{\alpha}{2} R^2t+2\alpha R^2\right],
\end{aligned}
\end{equation*}
where $V_f:=V_f(B_x(\sqrt{t}))$. Similar to the discussion of $J_1$,
by choosing $T$ sufficiently small, for all $t\in(0,T)$ and all $x\in M$,
$J_2$ also tends to zero when $R$ tends to infinity.\\

\emph{Step 4}. Recall that the co-area formula states that for all $h\in C^\infty_0(M)$,
\[
\int_{B_o(R+1)\backslash B_o(R)}h(y)d\mu(y)
=\int^{R+1}_R\left[\int_{\partial B_o(r)}h(y)d\mu_{\sigma,r}(y)\right]dr,
\]
where $\mu_{\sigma,r}$ denotes the weight area-measure induced by $\mu$
on $\partial B(o, r)$.
By the mean value theorem, for any $R>0$ there exists
$\bar{R}\in (R,R+1)$ such that
\begin{equation*}
\begin{aligned}
J:&=\int_{\partial B_o(\bar{R})}\left[H(x,y,t)|\nabla g|(y)+|\nabla H|(x,y,t)g(y)\right]d\mu_{\sigma,\bar{R}}(y)\\
&=\int_{B_o(R+1)\backslash B_p(R)}\left[H(x,y,t)|\nabla g|(y)+|\nabla H|(x,y,t)g(y)\right]d\mu(y)\\
&=J_1+J_2.
\end{aligned}
\end{equation*}
By step 2 and step 3, we know that by choosing $T$ sufficiently small,
 for all $t\in(o,T)$ and all $x\in M$,
$J$ tends to zero as $\bar{R}$ (and hence $R$) tends to infinity.
Therefore we complete Theorem \ref{liouL1} for $T$ sufficiently small.\\

\emph{Step 5}. At last, using the semigroup property of the $f$-heat
equation, we have
\begin{equation*}
\begin{aligned}
\frac{\partial}{\partial(s+t)}\left(e^{(s+t)\Delta_f}g\right)&=
\frac{\partial}{\partial t}\left(e^{s\Delta_f}e^{t\Delta_f}g\right)=
e^{s\Delta_f}\frac{\partial}{\partial t}\left(e^{t\Delta_f}g\right)\\
&=e^{s\Delta_f}e^{t\Delta_f}(\Delta_fg)=
e^{(s+t)\Delta_f}(\Delta_fg)
\end{aligned}
\end{equation*}
which implies Theorem \ref{liouL1} for all time $t$.
\end{proof}

Now we can finish the proof of Theorem \ref{liouL0}, following
the idea in \cite{[Li0]}.
\begin{proof}[Proof of Theorem \ref{liouL0}]
Let $g(x)$ be a nonnegative, $L^1$-integrable and $f$-subharmonic
function defined on $M$. Now we define
\[
g(x,t):=\int_MH(x,y,t)g(y)d\mu(y)
\]
with $g(x,0)=g(x)$. By Theorem \ref{liouL1},
\begin{equation*}
\begin{aligned}
\frac{\partial}{\partial t}g(x,t)
&=\int_M\frac{\partial}{\partial t}H(x,y,t)g(y)d\mu(y)\\
&=\int_M{\Delta_f}_yH(x,y,t)g(y)d\mu(y)\\
&=\int_MH(x,y,t){\Delta_f}_yg(y)d\mu(y)\geq 0.
\end{aligned}
\end{equation*}
Therefore we confirmed $g(x,t)$ is increasing for all $t$.
On the other hand, under the assumption of our theorem, by
Lemma 2.1 in \cite{[MuWa2]} we conclude that
\[
V_f(B_o(R))\leq C e^{c(n)R^2}
\]
holds for all $R>0$, where $C>0$ is a constant depending
on $A$ and $B_o(1)$. Hence
\[
\int^{\infty}_1\frac{R}{\log V_f(B_o(R))}dR=\infty.
\]
By Grigor'yan's result in \cite{[Grig2]} (see also Theorem 3.13 in \cite{[Grig3]}),
this implies
\[
\int_MH(x,y,t)d\mu(y)=1
\]
for all $y\in M$ and $t>0$. To finish our theorem, this equality implies
\[
\int_Mg(x,t)d\mu(x)=\int_M\int_MH(x,y,t)g(y)d\mu(y)d\mu(x)=\int_M g(y)d\mu(y).
\]
Since $g(x,t)$ is increasing in $t$, we conclude that $g(x,t)=g(x)$ and hence
$\Delta_fg(x)=0$. On the other hand, for any positive constant $a$,
let us define a new function
\[
h(x):=\min\{g(x),a\}.
\]
Then $h$ satisfies
\[
0\leq h(x)\leq g(x), \quad|\nabla h|\leq |\nabla g|\quad
\mathrm{and}\quad\Delta_fh(x)\leq 0.
\]
In particular, it will satisfy the same estimates, \eqref{mjbd} and \eqref{guji},
as $g$. Hence we can show that
\[
\frac{\partial}{\partial t}\int_MH(x,y,t)h(y)d\mu(y)
=\int_MH(x,y,t){\Delta_f}_yh(y)d\mu(y)\leq 0.
\]
Note that $h$ is still $L^1$, following the same argument as before,
we have $\Delta_fh(x)=0$. By the regularity theory of $f$-harmonic functions, this
is impossible unless $h=g$ or $h=a$. Since $a$ is arbitrary and $g$ is nonnegative,
this implies $g$ must be identically constant.
\end{proof}

\section*{Acknowledgments}
The work was initiated during the author visited the Academy of Mathematics
and System Sciences of the Chinese Academy of Sciences in Spring 2012.
The author would like to thank Professor X.-D. Li for his invitation
hospitality and helpful discussions. He brought this problem to the
author's attention. The author would also like to thank the anonymous
referee for pointing out a lot of expression errors and give many valuable
suggestions that helped to improve the presentation of the paper.

\bibliographystyle{amsplain}

\end{document}